\documentclass[11pt]{amsart}
\usepackage{graphics,psfrag}

\newcommand{\inv}{^{-1}} 
\newcommand{\tto}{\buildrel * \over\rightarrow }
\newcommand{\abs}[1]{\vert #1\vert} 
 
\newcommand{\edge}[1]{\buildrel#1 \over \rightarrow } 
\newcommand{\longedge}[1]{\buildrel#1 \over \longrightarrow } 
\newcommand{\set}[1]{\{ #1 \}}

\newcommand{\ovr}[1]{\overline #1}
\newtheorem{theorem}{Theorem}
\newtheorem{lemma}[theorem]{Lemma}

\theoremstyle{definition}

\theoremstyle{remark}

\title{On the Definition of Word Hyperbolic Groups}
\author{Robert H. Gilman}
\address{Department of Mathematical Sciences, Stevens Institute of 
Technology, Hoboken, New Jersey}
\email{rgilman@stevens-tech.edu}
\thanks{The author expresses his thanks to
the City College 
of New York for its hospitality while this paper was being written.}

\begin{document}

\begin{abstract}
Formal languages based on multiplication tables of finitely
generated groups are investigated and used to give a linguistic
characterization of word hyperbolic groups.
\end{abstract}

\maketitle

\section{Introduction}\label{intro}

Over the last several years combinatorial group theory has been
influenced by ideas from both low dimensional topology and formal
language theory. Applications of language theory include the
classification of groups with context--free word problem~\cite{MS1}
(together with~\cite{Du}), the use of indexed languages to describe
the fundamental groups of the known compact 3--manifolds~\cite{BG}, a
forthcoming complexity--theoretic analog of the Higman embedding
theorem~\cite{B+}, and the general theory of automatic
groups~\cite{E+}. In this paper we use formal languages in a novel way
to obtain a linguistic characterization of word hyperbolic groups.

A formal language is a subset of a free monoid $\Sigma^*$ over a
finite alphabet $\Sigma$. The connection between between a group $G$
and languages over $\Sigma$ is made by means of a surjective monoid
homomorphism $\Sigma^*\to G$ which maps $w\in \Sigma^*$ to
$\ovr w \in G$. The usual languages considered are the word problem,
$\{w\mid \ovr w=1\}$, and combings, i.e., languages projecting onto
$G$. We consider instead languages derived from the multiplication
table of $G$. For this purpose we need a new letter $\#$ not in the
alphabet $\Sigma$. 

\begin{theorem}\label{hyp} Let $\Sigma^*\to G$ be a choice of
generators for the group $G$. $G$ is word--hyperbolic if and only if
for some regular combing $R\subset \Sigma^*$, the language $M=\{u\#
v\# w\mid u,v,w \in R, \ovr u\ovr v\ovr w=1\}$ is context--free.
\end{theorem}

In short $G$ is hyperbolic if and only if it has a context--free
multiplication table. It is interesting that the original geometric
definition of word hyperbolic groups in terms of the thin triangle
condition is equivalent to a purely language-theoretic definition.
Choices of generators are defined in Section~\ref{defs}.

For any combing $R$ call $M=\{u\# v\# w\mid u,v,w \in R, \ovr u\ovr
v\ovr w=1\}$ the multiplication table determined by $R$. 
We investigate multiplication tables of virtually
free and automatic groups. Theorem~\ref{sigma} is a variation on the
main result of~\cite{MS1}.

\begin{theorem}\label{sigma} Let $\Sigma^*\to G$ be a choice of
generators and $M$ the multiplication table corresponding to the
combing $R=\Sigma^*$.
\begin{enumerate}
\item\label{finite} $G$ is finite if and only if $M$ is a regular
language.
\item\label{vf} $G$ is virtually free if and only if $M$ is
context--free.
\end{enumerate}
\end{theorem}

In Theorem~\ref{bi} we consider columns of the multiplication
table. The column of $g\in G$ is 
$C(g)=\{u\#w\mid u,w\in R, \ovr u g \ovr w = 1 \}$. Columns are
related to the comparator automata used in the definition of automatic
groups. Suppose $G$ is automatic with respect to the combing used to define
$M$, and $a\in \Sigma_\epsilon = \Sigma\cup\set\epsilon$ where
$\epsilon$ is the empty word. The binary relation
accepted by the comparator automaton for $a$ is $\{(u,w)\mid u,w\in R,
\ovr u\ovr a = \ovr w\}$ while $C(\ovr a)=\{u\#w \mid u,w\in R, \ovr
u\ovr a\ovr w=1\}$. 

\begin{theorem}\label{bi}Let $\Sigma^*\to G$ be a choice of generators. 
There exists a combing $R\subset \Sigma^*$ such that $C(\ovr a)$ is
context--free for all $a\in \Sigma_\epsilon$ if and only if $G$ is
asynchronously automatic with respect to a combing contained in $\Sigma^*$
and closed under taking formal inverses.
\end{theorem}

Groups which are asynchronously automatic with respect to a combing
closed under formal inverses form a subclass of asynchronously
biautomatic groups. It does not seem to be known whether or not this
subclass is proper.

The proof of Theorem~\ref{hyp} depends on the fact that the thin
triangle condition can be relaxed. The distance from a point on one
side of a triangle to the union of the other two sides may be allowed
to grow with the size of the triangle, and the sides of the triangle
need not be geodesics. See Theorem~\ref{combing} in
Section~\ref{flabby}.

Another linguistic characterization of hyperbolic groups is given by
Grunschlag~\cite[Section 3.2]{Gr}. He shows that hyperbolic groups are
those whose word problem is generated by a terminating growing
context--sensitive grammar.

Hyperbolic groups were introduced by Gromov~\cite{Gromov}. Additional
references are~\cite{A+}, \cite{Co} and \cite{Ghys}.

\section{Preliminary Items}\label{defs}

Keep the notation introduced in Section~\ref{intro}.
$G$ is a finitely generated group, $\Sigma$ is a finite alphabet, and
$\#$ is a letter not in $\Sigma$. $\Sigma_\# =
\Sigma\cup\set\#$, and $\Sigma_\epsilon=\Sigma\cup\{\epsilon\}$ where
$\epsilon$ stands for the empty word. $\Sigma$ has formal inverses if
it admits a permutation $a\to a\inv$ with orbits of length two. Formal
inverses on $\Sigma$ extend to formal inverses on $\Sigma^*$ by means
of the rule $(wv)\inv=v\inv w\inv$.

\subsection{Choice of Generators}\label{choice}

A choice of generators for $G$ consists of a finite alphabet $\Sigma$
equipped with formal inverses together with a surjective monoid
homomorphism $\Sigma^*\to G$ which maps $w\inv$ to $\ovr w \inv$.  Any
choice of generators $\Sigma^*\to G$ is extended to $\Sigma_\#^*\to G$
via $\ovr\#=1$.

Given a choice of generators $\Sigma^*\to G$ we define a path of
length $n$ in $G$ to be a sequence $g_0,\ldots, g_n$ of of group
elements such that $g_i=g_{i-1}\ovr{a_i}$ for some $a_i\in \Sigma$.
The label of this path is $a_1\cdots a_n\in \Sigma^*$. We use the
usual arrow notation $\cdots g_{i-1} \edge{a_i} g_i\cdots$. On
occasion we will allow paths with labels in $\Sigma^*_\#$.

Each word in $\Sigma^*_\#$ determines a path up to left--translation by
$G$. We identify words with paths and specify a particular path
corresponding to a word when necessary.  The
length of the shortest path from $g$ to $h$ is $d(g,h)$, a left--invariant
metric on $G$.  Shortest paths are called geodesics.

\subsection{Triangles}

A triangle $T$ consists of three points in $G$ joined by paths with
labels in $\Sigma^*$. These
paths are the sides of $T$. $T$ is $\delta$--thin if the distance from
any point on one side to the union of the other two sides is at most
$\delta$. The width of $T$, $\delta(T)$, is the smallest number
$\delta$ for which $T$ is $\delta$--thin. The norm of $T$, $\abs T$,
is the maximum distance between its vertices.  For any language
$L\subset\Sigma^*$, $T$ is an $L$-triangle if its sides are in $L$.
In particular $T$ is a geodesic triangle if its sides are
geodesics. If $\delta(T)$ is bounded as $T$ ranges over geodesic
triangles, then $G$ satisfies the thin triangle condition and is word
hyperbolic. 

Each element $u\#v\#w$ of the multiplication table corresponding to a
combing $R$ determines up to translation by $G$ an $R$--triangle with
sides $u,v,w$. For brevity we may refer to $u\#v\#w$ itself as a
triangle.

\subsection{Formal Languages}\label{formal}
See~\cite{Ha},\cite{HU},\cite{Re},\cite{Ro} for standard
introductions to the theory of automata and formal languages; a
group theoretic perspective is available in~\cite{Gilman2}.

Recall that context--free languages are the languages generated by
context--free grammars and that a context--free grammar $\mathcal G$
consists of a terminal alphabet $\Sigma$ (with or without formal
inverses), a set of nonterminals $N$, a 
start symbol $S\in N$, and a set of productions of the form
$A\to\alpha$ where $A\in N$ and $\alpha\in (\Sigma\cup N)^*$. All
these sets are finite.

Elements of $(\Sigma\cup N)^*$ are called sentential forms. The
notation $\alpha\to \beta$ means that the lefthand side of some
production is a subword of the sentential form $\alpha$, and that
the sentential form $\beta$ is obtained by replacing 
that subword by the righthand side of the production. The effect of
zero or more 
replacements is denoted by $\alpha \tto \beta$. When $\alpha\tto\beta$
we say $\beta$ is
derived from $\alpha$ or $\alpha$ derives $\beta$. The language generated
by $\mathcal G$ is $\set{w\mid w\in \Sigma^*, S\tto w}$.

\subsection{Transductions and Rational Subsets}

A rational transduction
$\rho:\Sigma^*\to \Delta^*$ from one finitely generated free monoid to
another is a rational subset of $\Sigma^*\times\Delta^*$. Write
$\rho(w)=v$ if $(w,v)\in \rho$ and $\rho(L)=\set{ v \mid \exists w\in
\Sigma^* \rho(w)=v}$ for $L\subset\Sigma^*$. The inverse of $\rho$ is
$\rho\inv=\set{(v,w)\mid (w,v)\in\rho}$.

The rational subsets of any monoid $P$ are the closure of its finite
subsets under union, product, and generation of
submonoid. Equivalently rational subsets are the subsets accepted by
finite automata over $P$. A finite automaton $\mathcal A$ over $P$ is
a finite directed graph with edge labels from $P$, a distinguished
initial state, and some distinguished terminal states. $\mathcal A$
accepts the set of labels of paths which begin at the initial state
and end at a terminal state.  Automata may be allowed to have more 
than one initial state. The accepted set is a union of sets accepted
by automata with unique initial states and so is rational.

Rational subsets of $\Sigma^*$ are called regular languages. It is
customary to restrict edge labels in automata over $\Sigma^*$ to
$\Sigma$ or $\Sigma_\epsilon$, but this restriction is not
necessary. Regular languages are closed under intersection and
difference while rational sets in general are not.

Since rational transductions are rational subsets, they are closed
under union, product and generation
of submonoids. They are also closed under inverse and under
composition in the sense of binary relations. Images of regular and
context--free languages under rational transductions are regular and
context--free respectively. In particular regular and context--free languages
are closed under homomorphism, inverse homomorphism, and
intersection with regular languages. Rational transductions are not
closed under intersection, but if $\rho:\Sigma^*\to \Delta^*$ is a
rational transduction, $R\subset \Sigma^*$ is regular, and $S\subset
\Delta^*$ is also regular, then $\rho\cap(R\times S)$ is a rational
transduction.   

\begin{lemma}\label{transduction} Fix $w,v\in
\Sigma^*$; $\rho=\set{(xwy, xvy) \mid x,y\in \Sigma^*}$ is a rational
transduction from $\Sigma^*$ to itself.
\end{lemma}
\begin{proof}
Let $D$ be the diagonal submonoid of
$\Sigma^*\times\Sigma^*$. Since $D$ is finitely generated, it is
rational. It follows that $\rho=D (w,v)D$ is a product of
rational sets and so is itself rational.
\end{proof}

\begin{lemma}\label{invert} If $\Sigma$ and $\Delta$ have formal
inverses and $\rho:\Sigma^*\to \Delta^*$ is a rational transduction,
then so is $\tau =\set{(w,v)\mid (w\inv, v\inv)\in \rho}$.
\end{lemma}
\begin{proof}
Pick an automaton accepting $\rho$. Reverse the orientation of each
edge and invert the edge label. Make every terminal state an initial
state and every initial state a terminal state.
\end{proof}

\begin{lemma}\label{linear} Let $\Sigma$ have formal inverses. A
relation $\rho:\Sigma^*\to \Sigma^*$ is a 
rational transduction if and only if $L=\set{ u\#w \mid (u, w\inv)\in
\rho}$ is generated by a context--free grammar with all productions of
the form $A\to xBy$ or $A\to x\#y$ for $x,y\in \Sigma^*$.
\end{lemma}
\begin{proof}
Suppose $\rho$ is accepted by a finite automaton $\mathcal A$ over
$\Sigma^*\times\Sigma^*$. Construct a context--free grammar with one
nonterminal $A_p$ for each vertex $p$ of $\mathcal A$. The start
symbol is the nonterminal corresponding to the initial vertex. For
each edge $p\longedge{(x,y)}q$ there is a production $A_p\to
xA_qy\inv$, and for each terminal vertex $q$ there is another
production $A_q\to \#$. It is straightforward to check that this
grammar generates $L$. The main step is
to use induction on path length and on derivation length to prove that
$A_p\tto uA_qw$ if and only if there is a path in $\mathcal A$ from
$p$ to $q$ with label $(u,w\inv)$. 

For the converse suppose $L$ is generated by a context--free grammar $\mathcal
G$ as above. Construct an automaton $\mathcal A$ whose vertices are the
nonterminals of $\mathcal G$ plus one terminal vertex. The initial
vertex is the start symbol. For each production $A\to xBy$ there is an
edge $A\longedge{(x,y\inv)} B$, and for each production  $A\to x\#y$ there
is an edge with label $(x,y\inv)$ from $A$ to the terminal vertex. Again
it is straightforward to check that $u\#w\in L$ if and only if
$\mathcal A$ accepts $(u,w\inv)$.
\end{proof}

\begin{lemma}\label{biauto}  $G$ is
asynchronously automatic with respect to a regular combing $R$ if
and only if for all $a\in 
\Sigma_\epsilon$ the relation $\rho_a = \set{(u,w)\mid u,w\in R,
\ovr u\ovr a=\ovr w}$ is a rational transduction.
\end{lemma}
\begin{proof} 
If $G$ is asynchronously automatic, then by~\cite[Definition
7.2.1]{E+} $\rho_a$ is a rational transduction. The converse is
~\cite[Theorem 1]{S} except that
the automata used there are more restricted than
ours. In terms of our notation the vertices of those automata are
partitioned into two sets. All edges leaving the first set have labels
from $\Sigma_\epsilon\times \set\epsilon$, and all edges leaving the
second set have labels from $\set\epsilon\times \Sigma_\epsilon$.

An automaton in our sense can be transformed into one satisfying the
definition in~\cite{S}.  Replace edges
by paths if necessary to insure that edge labels are from
$(\Sigma_\epsilon\times \set\epsilon) \cup (\set\epsilon\times
\Sigma_\epsilon) $. If vertex $p$ is a source for edges of both types,
add a vertex $p'$ and make all the edges of one type start at $p'$
instead of $p$. Add edges from $p$ to $p'$ and $p'$ to $p$ with label
$(\epsilon,\epsilon)$.

There is one more detail. In~\cite{S} $\rho_a$ is defined as
$\set{(u\$,w\$)\mid u,w\in R, \ovr{ua}=\ovr w }$ instead of
$\set{(u,w)\mid u,w\in R, \ovr{ua}=\ovr w }$, but if one version of
$\rho_a$ is a rational transduction, then the other one is too.
\end{proof}

\section{Flabby Triangles} \label{flabby}

In preparation for the proof of Theorem~\ref{hyp} we show that the
thin triangle condition used to define word hyperbolic
groups can be weakened.

\begin{theorem}\label{combing} A group $G$ is word--hyperbolic
if it admits a choice of generators $\Sigma^*\to G$ and a combing
$R\subset \Sigma^*$ such for some constant $C$ every $R$--triangle $T$
in $G$ has width $\delta(T)\le\abs T/75 +C$.
\end{theorem}

The rest of this section is devoted to proving
Theorem~\ref{combing}. Without loss of generality assume that there is
just one combing path for each $g\in G$. By~\cite[Theorem B]{Gilman}
it suffices to show that for some constant $K$ every cycle in $G$
of length $n$ can be triangulated with diagonals of length at most
$n/6 + K$. Before discussing triangulations we prove a lemma modeled
on~\cite[Lemma 1.5 of Chapter 3]{Co}.

\begin{lemma} \label{lemma}Let $w$ and $v_0$ be paths in
$G$ from $g$ to $h$ with $\abs w=n$ and $v_0\in R$.  There is a
constant $D$ independent of $n$ such that every point on $v_0$ is a
distance at most $n/36+D$ from $w$. If $w$ is a geodesic, then every
point of $w$ is a distance at most $n/18+2D$ from $v_0$.
\end{lemma}

\begin{proof}
Let $p_0$ be a point on $v_0$. If $n\le 2$, then the distance from
$p_0$ to $w$ is at most $E$, the maximum length of the finitely many
combing paths for elements $g\in G$ with $d(1,g)\le 2$.  Otherwise
estimate the distance by constructing a sequence of triangles as in
Figure~\ref{triangles}.
\begin{figure}[ht] 
\psfrag{g}[cc]{$g$} \psfrag{h}[cc]{$h$} \psfrag{w}[cl]{$w$}
\psfrag{v0}[bc]{$v_0$} \psfrag{v1}[cl]{$v_1$} \psfrag{p0}[cc]{$p_0$}
\psfrag{p1}[bc]{$p_1$} \psfrag{p2}[bc]{$p_2$}
\includegraphics{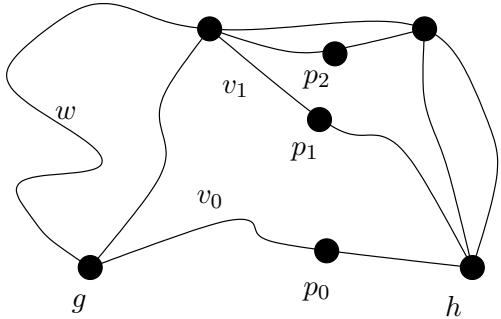}
\caption{A sequence of $R$--triangles.\label{triangles}}
\end{figure}
Let $T_0$ be an $R$--triangle whose base is $v_0$ and whose third
vertex is a point as close to the middle of $w$ as possible. One side,
call it $v_1$, of $T_0$ is distinct from $v_0$ and contains a point
$p_1$ with $d(p_1,p_0)\le\delta(T_0)$. If $v_1$ subtends a segment of
$w$ of length greater than $2$, construct triangle $T_1$ with base
$v_1$ in the same way $T_0$ was constructed. Continue until reaching a
triangle $T_m$ with point $p_{m+1}$ on a side subtending a segment of
$w$ of length at most $2$.

To show that the sequence of triangles terminates consider the
sequence of numbers defined by $b_0=n$ and $b_{k+1}=(1 +
b_k)/2$. From the construction above it is clear that the base of
triangle $T_k$ subtends a segment of $w$ of length at most $b_k$. It
is straightforward to show that $b_k\le 1+n/2^k$ whence the sequence
of triangles stops at $T_m$ for some $m\le\log_2(n)$.  

Since the third vertex of each $T_k$ lies on the segment of $w$
subtended by the base of $T_k$, we have $\abs{T_k}\le b_k$.
The distance from $p_0$ to $w$ is at most $E+\sum\delta(T_k) \le E +
\sum (b_k/75 +C)\le E+2n/75+(1+C)(1+\log_2 n)\le n/36 +D$ for some
constant $D$. 

To verify the last assertion of the lemma assume that $w$ is a geodesic and set
$A=n/36+D$. The points on $w$ a distance greater than $A$ from $v_0$
form a union of subpaths of $w$ not containing $g$ or $h$.  Let
$u$ be any such subpath, and write $w=w_1uw_2$. Observe that $w_1$
starts at $g$, $w_2$ ends at $h$, and each point of $v_0$ is
a distance at most $A$ from $w_1$ or $w_2$. It follows that
there are two adjacent points $g_1, g_2$ on $v_0$ and points $h_i$ on
the paths $w_i$ such that $d(g_i, h_i)\le A$. Since $w$ is a geodesic,
the distance along $w$ from $h_1$ to $h_2$ is $d(h_1,h_2) \le
d(h_1,g_1)+1+d(g_2,h_2)\le 2A +1$. But then any point on $u$ is a
distance at most $2A$ from $v_0$.
\end{proof}

We continue with the proof of Theorem~\ref{combing}.  Recall that
it suffices to triangulate $w$, a cycle of length $n$ in $\Gamma$,
with diagonals of length at most $n/6+K$. Take $K=6D+3$ where $D$ is
the constant from Lemma~\ref{lemma}.

To triangulate $w$ realize it as a regular polygon $P$ in the
plane. The vertices of $P$ are labelled by the group elements
$g_1,\ldots,g_n$ which occur along $w$, and the edges are labelled by
the letters of $w$. Particular group elements may occur more than once
as labels.  If $n\le 3$, then $w$ is deemed to be triangulated without
 any diagonals. Otherwise join the vertices of $P$ in pairs by
diagonals, i.e., straight line segments, so that no two diagonals meet
in the interior of $P$; the interior is divided into triangles; and
each edge of $P$ is one side of a triangle. The length of a diagonal
is the distance in $G$ between the labels of its endpoints. Edge
lengths are defined similarly and are either $0$ or $1$.

If $\abs w > 3$, $w$ can be triangulated in the following way so that
 all diagonals have length at most $n/6 + K$. First add a diagonal
 from $g_n$ to $g_2$; this diagonal has length at most 2.  We are done
 if $n=4$. Otherwise it suffices to show that whenever a diagonal of
 length at most $n/6 + K$ has endpoints $g_i, g_j$ with $3 \le j-i$,
 then we can add a diagonal from $g_i$ to $g_{j-1}$ or one from $g_{i+1}$
 to $g_{j}$ or diagonals from $g_i$ and $g_j$ to some $g_k$
 with $i+2 \le k \le j-2$. In other words it suffices that
 $d(g_i,g_k)$ and $d(g_j,g_k)$ are at both most $n/6 + K$ for some $k$
 with $i<k<j$.

Pick $h\in G$ as
close as possible to the middle of a geodesic path from $g_i$ to
$g_j$. By Lemma~\ref{lemma} $h$ is a distance at most $n/18+2D$ 
from some point on the $R$--path from $g_i$ to $g_j$, and that point
is itself a distance at most $n/36+D$ from the segment $w'$ of $w$
beginning at $g_i$ and ending at $g_j$. Consequently $h$ is a distance
at most $n/12+3D$ from $w'$, and it follows that $d(h,g_k)\le
n/12+3D+1$ for some $g_k$ with $i< k < j$. Hence $d(g_i,g_k)\le
d(g_i,h)+d(h,g_k)\le (1/2)(n/6+K +1) + n/12+3D+1\le n/6 +K$, and
likewise for $d(g_j,g_k)$.

\section{Hyperbolic Implies Context--Free}\label{hyperbolic}

With respect to any choice of generators $\Sigma\to G$ the geodesic
combing $R$ of a hyperbolic group $G$ is a regular
language~\cite[Theorem 13 in Chapter 9]{Ghys}. In this section we show
that the multiplication table $M$ determined by $R$ is context--free.

\begin{figure}[ht]
\psfrag{S}[Bl]{$S$} \psfrag{SS}[Bl]{$\# $} \psfrag{a1}[br]{$a_1$}
\psfrag{a2}[br]{$a_2$} \psfrag{ai}[br]{$a_i$}
\psfrag{aii}[br]{$a_{i+1}=a_{p-j+1}$} \psfrag{am}[Br]{$a_p$}
\psfrag{b1}[Bl]{$b_1$} \psfrag{bj}[Bl]{$b_j$}
\psfrag{bk}[Bl]{$b_{j+1}=b_{q-k+1}$} \psfrag{bn}[Bl]{$b_q$}
\psfrag{c1}[Bl]{$c_1$} \psfrag{ck}[Bl]{$c_k$}
\psfrag{cpp}[Bl]{$c_{r-1}$} \psfrag{cppp}[cl]{\parbox{2em}{\centering
\setlength{\baselineskip}{1ex}$c_{k+1}=c_{r-i+1}$}}
\psfrag{cp}[Bl]{$c_r$} \psfrag{g1}[bl]{$X_1$} \psfrag{g2}[Bl]{$X_2$}
\psfrag{gi}[Bl]{$X_i$} \psfrag{gp}[Bl]{$Y_j$} \psfrag{h1}[bc]{$Z_1$}
\psfrag{hk}[Bc]{$Z_k$} \psfrag{hkk}[Bl]{$Z_{k-1}$}
\includegraphics{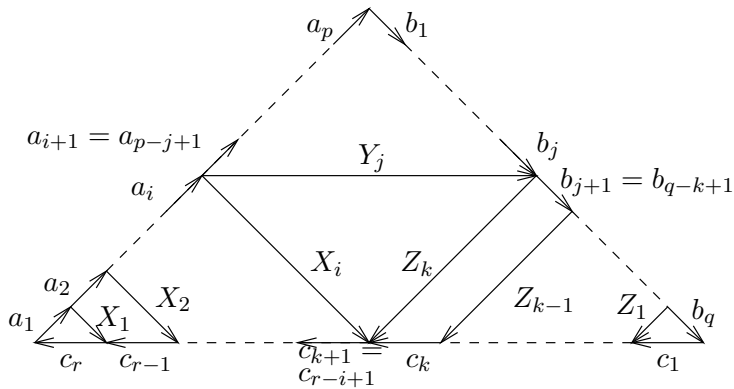}
\caption{A geodesic triangle.\label{triangle1}}
\end{figure}

Words $r\# s\# t\in M$ correspond to geodesic triangles whose sides
are paths with labels $r,s,t$.  By~\cite[Proposition 21 in Chapter
2]{Ghys} we may choose $\delta\ge 1$ so that points on the perimeter
of each geodesic triangle match in pairs with each point corresponding
to another point 
an equal distance from one of the vertices and matching points a
distance at most $\delta$ apart. Figure~\ref{triangle1} shows a
geodesic triangle with edge labels $a_1\cdots a_p$, The interior
arrows indicate paths in $G$ between matching points on the sides.  
The fact that points along the perimeter match in pairs
implies $i+j=p$, $j+k=q$, and $i+k=r$.
$b_1\cdots b_q$, $c_1\cdots c_r$. 

Define a context--free grammar $\mathcal G$ whose terminal alphabet is
$\Sigma_\#$ and whose nonterminal alphabet $N$ consists of a symbol
$X_w$ for each word $w\in \Sigma_\#^*$ of length $\abs w \le \delta$.
Let $V=\Sigma_\#\cup N$ and extend the choice of generators to a
monoid homomorphism $V^*\to G$ by $X_w\to \ovr{X_w}=\ovr w$ and
$\ovr\#=1$. The start symbol of $\mathcal G$ is $X_\epsilon$.

The productions of $\mathcal G$ are all replacements $X\to \alpha$
with $X\in V$, $\alpha$ a word of length at most 5 in $V$, and $\ovr X
= \ovr\alpha$. Since applying productions does not change images in
$G$, it is clear that $\mathcal G$ generates a context--free language
of words defining the identity in $G$.

To obtain a leftmost derivation of $a_1\cdots a_p\# b_1\cdots b_q\#
c_1\cdots c_r$ from Figure~\ref{triangle1} begin with $X_\epsilon\to
a_1X_1c_r$ and continue with productions corresponding to inscribed
quadrilaterals.
$$
X_\epsilon\to a_1X_1c_r\to a_1a_2 X_2c_{r-1}c_r \tto a_1\cdots
a_iX_ic_{r-i+1}\cdots c_r
$$
Apply the production $X_i\to Y_jZ_k$, and then do
$$
Y_j \to a_{i+1}Y_{j-1}b_{j}\tto a_{i+1}\cdots a_{p-1} Y_1
b_2\cdots b_{q-k} \to a_{i+1}\cdots a_p \# b_1\cdots b_{j}.
$$
Treat $Z_k$ similarly.

\begin{figure}[ht] 
\psfrag{a}{$a$} \psfrag{b}{$b$} \psfrag{c}{$c$} \psfrag{x}{$X$}
\psfrag{y}{$Y$} \psfrag{z}{$Z$} \includegraphics{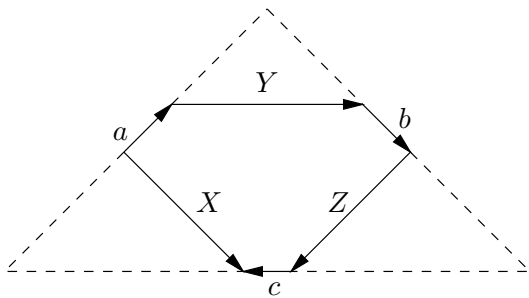}
\caption{Another geodesic triangle.\label{triangle2}}
\end{figure}

In the preceding derivation the right--hand sides of all productions
have length at 3. The reason we require productions with longer
right--hand sides is that in some geodesic triangles the central
figure is a hexagon instead of a triangle.  For that case we need
productions $X\to aYbZc$ as in
Figure~\ref{triangle2}. 

We see that $\mathcal G$ generates a context--free language $L$ which
contains $M$ and projects to $1$ in $G$. It follows that
$M=L\cap R\# R\# R$; and as intersections of context--free and regular
languages are context--free, $M$ is context--free.

\section{Context--Free Implies Hyperbolic}\label{cfsection}

In this section we complete the proof of Theorem~\ref{hyp}. Let
$\Sigma^*\to G$ be a choice of generators, $R$ a regular combing, and
$M$ the multiplication table determined by $R$. Assume
$M$ is context--free. By Theorem~\ref{combing} it suffices to show
that for some combing $R'$ every $R'$--triangle is $\delta$--thin.

Let $\mathcal G$ be a context--free grammar for $M$ in Chomsky normal
form. This normal form condition means that the  productions of
$\mathcal G$ look like $A\to BC$ 
or $A\to a$ where $a\in \Sigma_\#$ and $A,B,C$ are nonterminals. Without
loss of generality we may assume that each production participates in
some derivation of a word in $M$ and that each nonterminal occurs in a
production. For each nonterminal $A$ let $L_A$ be the context--free
language of all 
terminal words derived from $A$ by applying productions of $\mathcal
G$. Our conditions guarantee that $L_A$ is nonempty. Define $u_A$
to be a shortest word in $L_A$, and let $K$ be any constant greater
than the length of every $u_A$.

We claim that for a fixed nonterminal $A$ each word in $L_A$ represents
the same element of $G$. Indeed $A$ occurs in a derivation of some
$u\#v\#w\in M$ and derives a subword $x$ of $u\#v\#w$.  Because of the
way derivations are defined for context--free grammars, replacing $x$
by any $y\in L_A$ yields another member of $M$. As all elements of $M$
represent $1$ in $G$ (recall that $\ovr \#=1$), it follows that $\ovr x =
\ovr y$. The same 
reasoning shows that every word in $L_A$ contains the same number of
$\#$'s. Define that number to be the rank of $A$.

Fix an $R$--triangle $T$ with label $u\#v\#w\in M$;  $uvw$ is the
label of a cycle $g_1,\ldots,g_n$ in $G$, and subwords of $uvw$ are
paths.  We may think of each letter in $uvw$ as joining two group
elements in the cycle. 

\begin{figure}[ht]
\psfrag{u}[bc]{$u$} \psfrag{w}[bc]{$w$} \psfrag{b}[bc]{$b$}
\psfrag{ua}[bc]{$u_A$} \psfrag{wa}[bc]{$x_A$} \psfrag{wb}[bc]{$x_B$}
\psfrag{#}[bc]{$\#$}
\includegraphics{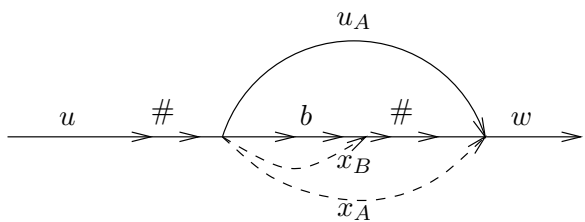}
\caption{Estimating $\delta$.\label{subwords}}
\end{figure}

Pick a letter $b$ in $uvw$. We will estimate the distance from the 
group elements it joins to another side of $T$. Any derivation of
$u\#v\#w$ can be written as $S\tto \alpha A\beta \to \alpha BC\beta
\tto u\#v\#w$ where $A$ is the 
last nonterminal of positive rank which derives a subword $x_A$ containing
$b$, and $B$ or $C$ is a nonterminal of rank zero deriving subword
of $x_A$ containing $b$. Assume $B$ has rank zero and derives a
subword (in fact a prefix) $x_B$ of $x_A$ containing $b$. The argument
is the same in the other case. The situation is illustrated in
Figure~\ref{subwords} for the case that $b$ lies in $v$. The dashed
lines in Figure~\ref{subwords} represent subwords of $u\#v\#w$. 

Since $x_A$ contains one or two $\#$'s, it begins on one side of $T$
and ends on another. There is a path $u_A$ with the same
initial and terminal point as $x_A$, and consequently the distance
from $b$ to another side of $T$ is at most $\abs{x_B} + \abs{u_A} \le
\abs{x_B} + K$. Thus the following lemma completes the proof
of Theorem~\ref{hyp}.

\begin{lemma}\label{thin} There exists a constant $K'$ and a regular combing
$R'\subset R$ such that every word in the restriction of $M$ to $R'$
affords a derivation in which each nonterminal of rank zero derives a
subword of length at most $K'$.
\end{lemma}

\begin{proof}
By the Pumping Lemma for context--free languages there exists a
constant  $K'$ such that
for all $A$ each $z \in L_A$ of length $\abs z\ge K'$ contains a
subword $x$ of length at most $K'$ which can be replaced by a shorter
word $y$ to obtain another element of $L_A$. Note that $\ovr x = \ovr
y$ as all words in $L_A$ represent the same element of $G$.

By Lemma~\ref{transduction}
$\rho_{x,y}=\set{(pxq, pyq)\mid p,q\in \Sigma_\#^*}$ is a rational
transduction. Hence so is $\rho$, the union of $\rho_{x,y}$ over all
pairs $x,y$ which occur above for nonterminals of rank 0.  Define
$R'=R-\rho\inv(R)$. $R'$ consists of all words in $R$ which cannot be
reduced to other words in $R$ by a substitution of the form $x\to
y$. Since these reductions are length reducing, $R'$ contains all
words in $R$ which are of minimal length among words in $R$ representing the
same element of $G$. Consequently $R'$ is still a combing. 
$R'$ is regular as it is a difference of regular languages.

Let $M'=M\cap (R'\#R'\#R')$, and consider a derivation
$S\tto u\#v\#w\in M'$. Each nonterminal of rank zero appearing in this
derivation derives a subword $z$ of $u$,$v$, or $w$. If $\abs z >
K'$, then one of $u,v,w$ can be
shortened by a substitution of $y$ for $x$ contrary to our choice of
$u\#v\#w$.  
\end{proof}

\section{The Combing $\Sigma^*$}

This section is devoted to the proof of Theorem~\ref{combing}.
Consider a choice of generators $\Sigma^*\to G$. Let $W=\set{w\mid
\ovr w = 1}$ be the word problem and $M$ be the multiplication table
determined by the combing $\Sigma^*$. It is well known that $G$ is
finite if and only if $W$ is regular (see~\cite{Gilman2} for example),
and by~\cite{MS1} together with~\cite{Du} $G$ is virtually free if and
only if $W$ is context--free. Thus it is enough to show $M$ is regular
if and only $W$ is and $M$ is context-free if and only if $W$ is.  We
give the argument for the regular case. The argument for the
context--free case is exactly the same.

Let $W_1$ be the inverse image of 1 under the extended homomorphism
$\Sigma_\#^* \to G$.  Observe that $W_1=f\inv(W)$ where
$f:\Sigma_\#^*$ to $\Sigma^*$ is defined by $f(\#)=\epsilon$ and
$f(a)=a$, $a\in\Sigma$.

Suppose $W$ is regular; then $M=W_1\cap\Sigma^*\#\Sigma^*\#\Sigma^*=
f\inv(W)\cap\Sigma^*\#\Sigma^*\#\Sigma^*$ is also regular.  Conversely
if $M$ is regular, then $W=f(W\#\#) = f(M \cap \Sigma^*\#\#)$ is
regular too.

\section{Automatic Groups}

In this section we prove Theorem~\ref{bi}. First suppose $G$ admits an
asynchronous automatic structure based on a combing $R$ which is
closed under taking inverses. For each $a\in\Sigma_\epsilon$ the
relation $\rho_a = \set{(u,v)\mid u,v\in R, \ovr u \ovr a = \ovr v}$ is
a rational transduction. By Lemma~\ref{linear} $L_a=\set{u\#w\mid
u,w\inv \in R, \ovr u \ovr a = \ovr w\inv }$ is context--free. But
as $R$ is closed under inverses, $L_a=\set{u\#w\mid
u,w \in R, \ovr u \ovr a \ovr w = 1 } = C(\ovr a)$.

To prove the converse fix $a\in \Sigma^*_\epsilon$, 
pick a context--free grammar $\mathcal G$ for $C(\ovr a)$ in Chomsky
normal form, and argue as in
Section~\ref{cfsection}. We may assume that $\mathcal G$ has no superfluous
nonterminals or productions and consequently that the nonterminals of
$\mathcal G$ each have a well defined rank of zero or one. By expanding
nonterminals of rank one first we can put every derivation into the form  
$S\tto \alpha B \beta \tto u\# w$
where $S$ and $B$ have rank one, $\alpha$ and $\beta$ are words in
nonterminals of rank zero, $\alpha\tto u$, $\beta\tto w$, and $B\to
\#$. Further for some regular subcombing $R_a\subset R$, $u\#
w\in C(\ovr a)\cap (R_a\times R_a)$ implies that subwords of $u$ or
$w$ derived from nonterminals in $\alpha$ or $\beta$
have length at most $K$. 

Consider the linear context--free grammar $\mathcal G'$ obtained by
replacing each production $A\to BC$ of $\mathcal G$ in which $A$ and
$C$ are of rank one and $B$ of rank zero by the productions $A\to xC$
where $x$ ranges over all words of length at most $K$ in
$L_B$. Likewise productions $A\to BC$ with $A$ and $B$ of rank
one and $C$ of rank zero are replaced by productions $A\to By$, $y\in
L_C$, $\abs y \le K$. Let $L$ be the language generated  by $\mathcal
G'$. Clearly $L\subset C(\ovr a)$, and from the discussion above it
follows that $ C(\ovr a)\cap (R_a\times R_a)\subset L$.

Note that nonterminals of rank zero do not appear in any $\mathcal
G'$--derivations of words in $L$.  After all nonterminals of rank zero
are deleted $\mathcal 
G'$ satisfies the hypothesis of Lemma~\ref{linear}, and it follows
that $\tau_a=\set{(u,v)\mid u\#v\inv\in L}$ is a rational
transduction.

By construction each $R_a$ contains all words of minimal length among
those in $R$ defining the same element of $G$. Thus $R_1=\cap R_a$ is
a regular combing, and so is $R'=R\cup R_1\inv$. Replace each
$\tau_a$ by its restriction to $R_1\times R_1\inv$; $\tau_a$ is still
a regular transduction.  Let $\mu=\tau_\epsilon\inv$, and check that
$\tau_a \cup (\tau_a\circ\mu) \cup (\mu\circ\tau_a) \cup
(\mu\circ\tau_a\circ\mu) = \set{(u,v)\mid u,v\in R', \ovr u \ovr a =
\ovr v}$. It follows from Lemma~\ref{biauto} that the combing $R'$ supports
an asynchronous automatic structure for $G$.

\end{document}